\documentclass[11pt]{article}
\usepackage{latexsym}
\usepackage{a4,epsfig}
\usepackage{graphicx}
\usepackage{amscd}
\usepackage{amsbsy}
\newcommand\eq[1] {(\ref{#1})}

\newcommand{\nonum}{\nonumber \\}
\newcommand{\beqa}{\begin{eqnarray}}
\newcommand{\eeqa}[1]{\label{#1}\end{eqnarray}}
\newcommand{\beq}{\begin{equation}}
\newcommand{\eeq}[1]{\label{#1}\end{equation}}
\newcommand{\R}{\mbox{R}}

\newcommand{\Tr}{\mathop{\rm Tr}\nolimits}

\newcommand{\Md}{\partial}

\newcommand{\qed}{\hfill \ensuremath{\Box}}
\newcommand{\GO}{\Omega}

\newcommand{\CB}{{\mathcal B}}

\newcommand{\CS}{{\mathcal S}}

\newcommand{\ds}{\displaystyle}

\def \ba {\begin{array}}
\def \ea {\end{array}}
\def \refe #1.{(\ref{#1})}
\def \reff #1.{figure~\ref{#1}}
\def \refs #1.{section~\ref{#1}}
\def \refss #1.{subsection~\ref{#1}}
\def \refD #1.{Definition~\ref{#1}}
\def \refT #1.{Theorem~\ref{#1}}
\def \refL #1.{Lemma~\ref{#1}}
\def \refC #1.{Corollary~\ref{#1}}
\def \refP #1.{Proposition~\ref{#1}}
\def \refR #1.{Remark~\ref{#1}}
\def \refE #1.{Example~\ref{#1}}
\def \refN #1.{Notation~\ref{#1}}

\usepackage{graphics}
\begin{document}
\vspace{-1in}
\title{Bounds on the volume fraction of 2-phase, 2-dimensional elastic bodies and on (stress, strain) pairs in composites}
\author{Graeme Walter Milton and Nguyen Hoang Loc\\
\small{Department of Mathematics, University of Utah, Salt Lake City UT 84112, USA} \\
\small{(milton@math.utah.edu, loc@math.utah.edu)}}
\date{}
\maketitle
\begin{abstract}
Bounds are obtained on the volume fraction in a two-dimensional body containing two elastically isotropic materials with known bulk and shear
moduli. These bounds use information about the average stress and strain fields, energy, determinant of the stress, 
and determinant of the displacement gradient, which can be determined from measurements
of the traction and displacement at the boundary. The bounds are sharp if in each phase certain displacement field components are constant. The inequalities we obtain
also directly give bounds on the possible (average stress, average strain) pairs in a two-phase, two-dimensional, periodic or statistically homogeneous
composite.
\end{abstract}
\vskip2mm

\noindent
\vskip2mm

%%%%%%%%%%%%%%%%%%%%%%%%%%%%%%%%%%%%%%%%

%ENVIRONMENTS THEOREMS...
% These are predefined, and follow the format of the journal!
%English
\newtheorem{theorem}{Theorem}[section]
\newtheorem{lemma}[theorem]{Lemma}
\newtheorem{e-proposition}[theorem]{Proposition}
\newtheorem{corollary}[theorem]{Corollary}
\newtheorem{e-definition}[theorem]{Definition\rm}
\newtheorem{remark}{\it Remark\/}
\newtheorem{example}{\it Example\/}
%%%%%%%%%%%%%%%%%%%%%%%%%%%%%%%%%%%%%%%%%%%%%%%%%%%%%%%%%%%%%%%%%%%%%%%%%%
\section{Introduction}
\label{intro}
\def\theequation{1.\arabic{equation}}\makeatother
\setcounter{equation}{0}
%%%%%%%%%%%%%%%%%%%%%%%%%%%%%%%%%%%%%%%%%%%%%%%%%%%%%%%%%%%%%%%%%%%%%%%%%%%%
A fundamental problem is to determine the volume fraction occupied by an inclusion in a body, or more generally
the volume fraction occupied one material in a body containing two materials. This can usually be done by weighing the body
but this may not always be practical or the densities of the two materials may be close. Then one might
seek to bound the volume fraction
from measurements of tractions and displacements
(for elasticity) or current fluxes and voltages (for conductivity)  at the boundary of the body. If the body contains a statistically homogeneous or periodic composite (with microstructure
much smaller than the dimensions of the body), then such boundary measurements can yield information about the effective tensors of the composite and it has long been recognized
(see, for example, \cite{Phan-Thien:1982:PUB,McPhedran:1982:ESI,McPhedran:1990:ITP,Cherkaeva:1998:IBM}) that bounds on effective tensors (which involve the volume fraction and material moduli) can be inverted to yield bounds on volume fractions. 
As shown in \cite{Milton:2011:UBE}, even if the body does not contain a 
statistically homogeneous or periodic composite, but provided that the applied tractions (or boundary displacements), or current fluxes (or boundary voltages) are
such that the fields in the body would be uniform were it filled with a homogeneous material, then boundary measurements can yield information about the effective tensor
of a composite containing rescaled copies of the body packed to fill all space in a periodic structure. Bounds on this effective tensor yield universal bounds on the
response of the body when such special boundary conditions are applied, which generalize those first obtained by Nemat-Nasser and Hori 
\cite{Nemat-Nasser:1993:MOP,Nemat-Nasser:1995:UBO}. They can then be inverted to yield
bounds of the volume fraction \cite{Milton:2011:UBE}, and when the volume fraction is asymptotically small the resulting bounds include those obtained by Capdeboscq and Vogelius 
\cite{Capdeboscq:2003:OAE,Capdeboscq:2004:RSR} (for conductivity)
and Capdeboscq and Kang \cite{Capdeboscq:2008:IHS} (for elasticity) using polarizability tensor bounds \cite{Lipton:1993:IEE,Capdeboscq:2003:OAE}.
Other bounds on the volume fraction, involving constants which are not easy to determine, were obtained by Kang, Seo and Sheen \cite{Kang:1997:ICP},
Ikehata \cite{Ikehata:1998:SEI}, Alessandrini and Rosset \cite{Alessandrini:1998:ICP}, Alessandrini, Rosset and Seo \cite{Alessandrini:2000:OSE},
and Alessandrini, Morassi and Rosset \cite{Alessandrini:2002:DCE}.

Given this close connection between bounding the effective tensor of a composite material and bounding the response of the body when these special boundary conditions are used, one
might wonder if methods that are used to obtain bounds on effective tensors of composites could also be used to bound the response of the body with any boundary conditions on the fields,
and then be inverted to bound the volume fraction. For conductivity such an approach has been successfully taken by Kang, Kim, and one of the authors \cite{Kang:2011:SBV}, using the translation method of
Murat and Tartar \cite{Tartar:1979:ECH,Murat:1985:CVH,Tartar:1985:EFC} and Lurie and Cherkaev \cite{Lurie:1982:AEC,Lurie:1984:EEC} 
which is one of the most successful methods for bounding effective tensors of composites: see the books \cite{Cherkaev:2000:VMS,Allaire:2002:SOH,Milton:2002:TOC,Tartar:2009:GTH}.
For a certain class of inclusion shapes (such that the field inside the inclusion is uniform for appropriate boundary conditions) one of the resulting volume fraction bounds gives the exact
volume fraction. For special boundary conditions the bounds reduce to those obtained in \cite{Milton:2011:UBE}, and for asymptotically small volume fractions the bounds reduce to those of 
Capdeboscq and Vogelius \cite{Capdeboscq:2003:OAE,Capdeboscq:2004:RSR}.  
 
The initial goal of this paper was to use the translation method to extend \cite{Kang:2011:SBV} to elasticity, using measurements obtained under a single, but arbitrary, loading,
and that is essentially done in Section \ref{A bound obtained by using the null-Lagragian}. The attainability conditions for the resulting bounds obtained in
Section \ref{The attainability condition for 2.15} then lead us to a new method for obtaining bounds which is not based on variational principles. This method, the method of splitting,
is described in Section \ref{The method of splitting}. It correlates volume averages of various quantities over each phase and then projects the information onto
the quantities of interest to obtain the desired bounds. This approach is likely to have wider applications, and in particular could lead to new bounds on the response of
(possibly non-linear) multiphase bodies for many properties, not just for elasticity.

The bounds we derive also directly give bounds on the possible (average stress, average strain) pairs in a two-phase, two-dimensional, composite. These bounds are the natural
generalization to elasticity of the conductivity bounds on possible (average electric field, average current field) pairs obtained by Ra{\u\i}tum 
\cite{Raitum:1978:EEP,Raitum:1983:QES} and Tartar \cite{Tartar:1995:RHM} (see also chapter 22.4 in \cite{Milton:2002:TOC}), which can also be generalized 
to non-linear materials \cite{Milton:2000:BCN,Talbot:2004:BEC,Peigney:2005:PBM,Peigney:2008:RSC,Bhattacharya:2010:SIP}. Only when one of the
phases is void, has the complete characterization of  possible (average stress, average strain) pairs been obtained \cite{Milton:2003:RAS}. In principle bounds
on the possible (average stress, average strain) pairs could be obtained from knowledge of the G-closure of all possible effective elasticity tensors associated
with composites of the two phases mixed in prescribed proportions. However, this G-closure is only partly known 
(for a survey of results see \cite{Cherkaev:2000:VMS,Allaire:2002:SOH,Milton:2002:TOC}.)

%%%%%%%%%%%%%%%%%%%%%%%%%%%%%%%%%%%%%%%%%%%%%%%%%%%%%%%%%%%%%%%%%%%%%%%%%%
\section{Preliminaries}
\label{prelim}
\def\theequation{2.\arabic{equation}}\makeatother
\setcounter{equation}{0}
%%%%%%%%%%%%%%%%%%%%%%%%%%%%%%%%%%%%%%%%%%%%%%%%%%%%%%%%%%%%%%%%%%%%%%%%%%%%

Let $\Omega$ be a smooth bounded domain in $\R^2,$ occupied by a two dimensional elastic body.
Assume that the body is made from two different isotropic elastic materials, characterized by their bulk moduli $\kappa_1$ and $\kappa_2$ and shear moduli $\mu_1$ and $\mu_2$. 
Denote, respectively, by $\sigma$ and $\varepsilon$ the stress and strain fields acting on the body under consideration. These fields are governed by the following equations.
\beq
	\varepsilon = S \sigma=(1/2\mu)\sigma+(1/4\kappa-1/4\mu)(\Tr\sigma)I, \quad \nabla \cdot \sigma = 0, \quad \varepsilon = \frac{1}{2}(\nabla u + \nabla u^T),
\eeq{1.1} 
where $u$ is the displacement field, $\CS$ is the compliance tensor, $I$ is the second order identity tensor, $\kappa$ is the bulk modulus (taking values $\kappa_1$ or $\kappa_2$)
and $\mu$ is the shear modulus (taking values $\mu_1$ or $\mu_2$). Although two dimensional bodies do not occur in practice, this formulation is applicable to problems of
plane stress or plane strain. 

It is assumed that one can measure the traction $\sigma\cdot n$ and the displacement $u$ at the boundary of $\GO$, where $n$ is the unit outward normal vector. From these measurements one can determine 
the volume averages of certain quantities. These quantities are null-Lagrangians, functionals of $u$ and/or $\sigma$ which can be integrated by parts, and expressed only in terms of boundary values: thus
their Lagrangian vanishes. Introducing angular brackets to denote a volume average, i.e. 
\beq \langle g \rangle= \frac{1}{|\GO|}\int_\GO g, \eeq{}
and choosing notations so that $\nabla u$ has $\nabla u_1$ and $\nabla u_2$ as its first and second columns, 
rather than rows, the five  null-Lagrangians we will work with are the average fields,
\beq \langle \sigma \rangle= \frac{1}{|\GO|}\int_{\Md\GO}x(\sigma\cdot n)^T, \quad
\langle \nabla u \rangle= \frac{1}{|\GO|}\int_{\Md\GO}n(u)^T, \quad
\eeq{}
the energy,
\beq \langle \sigma\cdot \varepsilon \rangle= \frac{1}{|\GO|}\int_{\Md\GO}(\sigma\cdot n)\cdot u,
\eeq{}
and the two additional null-Lagrangians
\beq a=\langle \det\sigma\rangle, \quad b=\langle \det \nabla u \rangle.
\eeq{}
To express the last two quantities in terms of boundary values, it is helpful to let $j_1$ and $j_2$ denote the
divergence free vector fields which are the first and second columns of $\sigma$, and to 
introduce the matrix
\beq R_\perp=\pmatrix{0 & 1 \cr -1 & 0} \eeq{}
for a clockwise $90^\circ$ rotation. Then, as follows directly from the analysis in \cite{Kang:2011:SBV},
\beq a=\langle j_1\cdot R_\perp j_2\rangle=\frac{1}{|\GO|}\int_{\Md\GO} q_1(x)(\int_{x_0}^x q_2),
\eeq{1.1a}
where $q_1$ and $q_2$ are the fluxes $q_1=j_1\cdot n$ and $q_2=j_2\cdot n$, which are components of the
traction $\sigma\cdot n$, $x_0\in\Md\GO$ and the last integral in \eq{1.1a} is along the boundary
$\Md\GO$ in the counterclockwise direction. Also, as follows directly from the analysis in \cite{Kang:2011:SBV},
\beq b=\langle \nabla u_1\cdot R_\perp \nabla u_2\rangle
=\frac{1}{|\GO|}\int_{\Md\GO}u_1 n\cdot R_\perp \nabla u_2=\frac{1}{|\GO|}\int_{\Md\GO}u_1\frac{\Md u_2}{\Md t},
\eeq{}
where $\Md/\Md t$ denotes the tangential derivative along $\Md\GO$ in the counterclockwise direction.

These are not the only null-Lagrangians. For $i=1, 2$ 
\beq  \langle \sigma\cdot \nabla u_i \rangle= \frac{1}{|\GO|}\int_{\Md\GO}(\sigma\cdot n)u_i
\eeq{}
is also a null-Lagrangian, but to simplify the analysis we refrain from considering these null-Lagrangians 
except in that linear combination which gives the energy.

Our objective is to find inequalities (bounds) which link the values of the five null-Lagrangians 
with the volume fraction in the body and the moduli of the materials. When the body is the unit
cell of a periodic composite material and periodic boundary conditions on the fields are imposed 
then the values of the average stress $\langle \sigma \rangle$ and average displacement gradient
$\langle \nabla u \rangle$ determine the values of the energy and $a$ and $b$:
\beq \langle \sigma\cdot \varepsilon \rangle=\langle \sigma\rangle\cdot\langle\varepsilon\rangle,
\quad a=\det \langle\sigma\rangle, \quad b=\det\langle\nabla u \rangle,
\eeq{comp}
as can be shown using Fourier analysis (see, for example, section 13.3 in \cite{Milton:2002:TOC}).
Also without loss of generality, by making an infinitesimal global rotation if necessary, we can 
assume that $\langle \nabla u \rangle$ is symmetric, in which case it can be identified with the
average strain $\langle\varepsilon\rangle$. Thus for composites the inequalities we obtain 
give bounds on the possible (average stress, average strain) pairs, and only incorporate the
volume fractions of the two isotropic phases and their moduli.

In Cartesian coordinates, $\sigma$ and $\varepsilon$ can be represented as  $2 \times 2$ matrices. It is convenient to use the basis
\beq
	\CB = \left\{
			\frac{1}{\sqrt{2}}\left[
								\begin{array}{cc}
									1&0\\0&1
								\end{array}
							  \right],
			\frac{1}{\sqrt{2}}\left[
								\begin{array}{cc}
									1&0\\0&-1
								\end{array}
							  \right],
			\frac{1}{\sqrt{2}}\left[
								\begin{array}{cc}
									0&1\\1&0
								\end{array}
							  \right]
		  \right\},
\eeq{1.3}
so that an arbitrary symmetric matrix 
$
	A = \left[
			\begin{array}{cc}
				a_{11}&a_{12}\\a_{12}&a_{22}
			\end{array}
		\right]
$ is represented by $v = \frac{1}{\sqrt{2}}(a_{11} + a_{22}, a_{11} - a_{22}, 2a_{12}).$ Thus, from now on we understand $\sigma$ and $\varepsilon$ as 3 dimensional vectors. 
If a $3-$dimensional vector $v$ represents a $2 \times 2$ symmetric matrix $A$ in this basis then its determinant is given by
\beq
	\det A = \frac{1}{2}v \cdot T v,
\eeq{det v} where
\beq
	T = \left[
		\begin{array}{ccc}
			1&0&0\\
			0&-1&0\\
			0&0&-1
		\end{array}
	\right].
\eeq{}
The null Lagrangian $\det \sigma$ plays a pivotal role in developing the bounds of the next section.

Recall that $\kappa$ and $\mu$ denote the bulk and shear, respectively, moduli of the elastic body under consideration. Then, in the basis $\CB$, 
the tensor $S$ can be expressed by the matrix
\beq
	S = \frac{1}{2}\left[
				\begin{array}{ccc}
					\frac{1}{\kappa}&0&0\\
					0&\frac{1}{\mu}&0\\
					0&0&\frac{1}{\mu}
				\end{array}
			\right],
\eeq{1.2}
because the elastic body is isotropic. 
%%%%%%%%%%%%%%%%%%%%%%%%%%%%%%%%%%%%%%%%%%%%%%%%%%%%%%%%%%%%%%
\section{Bounds obtained by the translation method using the null-Lagragian $\det \sigma$} \label{A bound obtained by using the null-Lagragian}
\def\theequation{3.\arabic{equation}}\makeatother
\setcounter{equation}{0}
%%%%%%%%%%%%%%%%%%%%%%%%%%%%%%%%%%%%%%%%%%%%%%%%%%%%%%%%%%%%%%%

Assume, in this section, that $\langle \sigma \cdot \varepsilon \rangle,$  $\sigma_0 = \langle \sigma \rangle,$ $\varepsilon_0 = \langle \varepsilon \rangle$ and $a = \langle \det \sigma \rangle$ can be evaluated or estimated. 
For $i = 1, 2$, let $\kappa_i$ and $\mu_i$ be the bulk and shear, respectively, moduli of the $i-$th phase. Then,
\beq
	 \kappa = \chi_1\kappa_1 + \chi_2\kappa_2 \mbox{ and } \mu = \chi_1 \mu_1 + \chi_2\mu_2,
\eeq{} where 
\beq
	\chi_i = \left\{
					\begin{array}{ll}
						1 & \mbox{ in phase } i,\\
						0 & \mbox{ otherwise}.
					\end{array}
				\right.
\eeq{}

Define
\beq
	\mu^* = \max\{\mu_1, \mu_2\}, \quad \kappa^* = \max\{\kappa_1, \kappa_2\},
\eeq{}
and fix $\alpha \in (-\frac{1}{2\mu^*}, \frac{1}{2\kappa^*})$ to ensure that the translated tensor
\beq
	L = S - \alpha T
\eeq{2.4}
is positive definite. The classical complementary energy minimization principle implies
\beq
	\langle\sigma \cdot \varepsilon \rangle = \min_{
															\begin{array}{c}
																\langle\underline\sigma\rangle = \sigma_0\\
																\langle S\underline\sigma\rangle = \varepsilon_0\\	
																\nabla \cdot \underline\sigma = 0\\
																\langle\det\underline\sigma\rangle = a \\
                                                                                                                                \underline\sigma\cdot n=\sigma\cdot n~~{\rm on}~~\Md\GO
															\end{array}
													   }
	\langle\underline \sigma \cdot S \underline \sigma\rangle,
\eeq{2.4a}
where the additional constraints that $\langle\underline\sigma\rangle = \sigma_0$, $\langle S\underline\sigma\rangle = \varepsilon_0$, 
and $\langle\det\underline\sigma\rangle = a$ have been added since we know this information about the minimizing fields. In the third constraint in \eq{2.4a} the divergence of $\underline \sigma$ 
is understood as the divergence of the matrix that $\underline \sigma$ represents. Using \eq{det v} we then have
\beq
	\langle\sigma \cdot \varepsilon \rangle - 2\alpha a = \min_{
															\begin{array}{c}
																\langle\underline\sigma\rangle = \sigma_0\\
																\langle S\underline\sigma\rangle = \varepsilon_0\\	
																\nabla \cdot \underline\sigma = 0\\
																\langle\det\underline\sigma\rangle = a\\
                                                                                                                 \underline\sigma\cdot n=\sigma\cdot n~~{\rm on}~~\Md\GO
															\end{array}
													   }
	\langle\underline \sigma \cdot L\underline \sigma\rangle.
\eeq{2.5}
Dropping the last three constraints in the minimum above and defining
\beq
	e_0 = \varepsilon_0 - \alpha T\sigma_0,
\eeq{2.6}
gives
\beq
	\langle\sigma \cdot \varepsilon \rangle - 2\alpha a \geq \min_{
															\begin{array}{c}
																\langle\underline\sigma\rangle = \sigma_0\\
																\langle L \underline\sigma\rangle = e_0
															\end{array}
													   }
	\langle\underline \sigma \cdot L\underline \sigma\rangle.
\eeq{30}

This minimum can be found using the Lagrange multiplier method. In fact, if $\hat \sigma$ denotes the minimizer of the right-hand side of \eq{30}, then there exist two constant vectors $\lambda_1$ and $\lambda_2$ such that
\beq
	2L \hat \sigma = \lambda_1 + L \lambda_2.
\eeq{2.7}
Denoting by $\langle \sigma \rangle_1$ and $\langle \sigma \rangle_2$ the average of $\sigma$ on phase $1$ and phase $2$, respectively, from \eq{2.7}, we have   
\beqa
			2L_1 \langle \sigma \rangle_1 &=& \lambda_1 + L_1\lambda_2,\nonum
			2L_2 \langle \sigma \rangle_2 &=& \lambda_1 + L_2\lambda_2,
\eeqa{}
which gives
\beqa
			\lambda_1 &=& 2(L_1 - L_2)^{-1}L_1L_2(\langle \sigma \rangle_2 - \langle \sigma \rangle_1),\nonum
			\lambda_2 &=& 2(L_1 - L_2)^{-1}(L_1\langle \sigma \rangle_1 - L_2\langle \sigma \rangle_2).
\eeqa{2.10}
On the other hand, we can see that
\beqa
			\langle \sigma \rangle_1 &=& \frac{1}{f_1}(L_1 - L_2)^{-1}(e_0 - L_2\sigma_0),\nonum
			\langle \sigma \rangle_2 &=& -\frac{1}{f_2}(L_1 - L_2)^{-1}(e_0 - L_1\sigma_0),
\eeqa{2.11} by solving the system
\beqa
			\sigma_0 &=& f_1\langle \sigma \rangle_1 + f_2\langle \sigma \rangle_2\nonum
			e_0 &=& f_1L_1\langle \sigma \rangle_1 + f_2L_2\langle \sigma \rangle_2.
\eeqa{}
Thus, we have
\beqa
			\lambda_1 &=& \frac{2}{f_1f_2}(L_1 - L_2)^{-2}L_1L_2(e_0 - \langle L \rangle \sigma_0),\nonum
			\lambda_2 &=& \frac{2}{f_1f_2}(L_1 - L_2)^{-2}[(f_2L_1 + f_1L_2)e_0 - L_1L_2\sigma_0].	
\eeqa{}

For simplicity in calculations, assume first that $\alpha \in (-\frac{1}{2\mu^*}, \frac{1}{2\kappa^*})$ so that $L$ is invertible. Later we will consider the
limiting case where $L$ is singular, but still positive semi-definite. It follows from \eq{2.7} that  
\beqa
			2 e_0 &=& \lambda_1 + \langle L \rangle \lambda_2,\nonum
			2 \sigma_0 &=& \langle L^{-1}\rangle \lambda_1 + \lambda_2.
\eeqa{2.8}
Moreover, we have
\beqa
	\langle \hat \sigma \cdot L \hat \sigma \rangle &=& \frac{1}{4}\langle (L^{-1} \lambda_1 +  \lambda_2) \cdot ( \lambda_1 + L \lambda_2) \rangle \nonum
	&=&\frac{1}{4} [\langle L^{-1} \rangle  \lambda_1 \cdot  \lambda_1 + 2  \lambda_1 \cdot  \lambda_2 + \langle L \rangle  \lambda_2 \cdot  \lambda_2].
\eeqa{}
On the other hand, it follows from \eq{2.8} that
\beq
	 \lambda_1 \cdot  \lambda_2 = 2 s_0 \cdot  \lambda_2 - \langle L \rangle  \lambda_2 \cdot  \lambda_2 = 2 \sigma_0 \cdot  \lambda_1 - \langle L^{-1} \rangle  \lambda_1 \cdot  \lambda_1.
\eeq{}
Thus, it follows that
\beqa
	\langle \hat \sigma \cdot L \hat \sigma \rangle &=& \frac{1}{2}[ \sigma_0 \cdot  \lambda_1 +  \lambda_2 \cdot  e_0]\nonum
	&=& \frac{1}{2}[(2 e_0 - \langle L \rangle \lambda_2) \cdot \sigma_0 + \lambda_2 \cdot e_0]\nonum
	&=& e_0 \cdot \sigma_0 + \frac{(L_1 - L_2)^{-2}}{f_1f_2}[e_0 - \langle L \rangle \sigma_0]\cdot[(f_2L_1 + f_1L_2)e_0 - L_1L_2\sigma_0].
\eeqa{}
This, together with \eq{2.6} and \eq{30}, gives us the bound
\beq
	\langle\sigma \cdot \varepsilon \rangle - \sigma_0 \cdot \varepsilon_0 - 2\alpha a + 2\alpha \det \sigma_0 \geq \frac{(L_1 - L_2)^{-2}}{f_1f_2}[e_0 - \langle L \rangle \sigma_0]\cdot[(f_2L_1 + f_1L_2)e_0 - L_1L_2\sigma_0].
\eeq{2.12}
By taking limits, we can see that \eq{2.12} is valid for $\alpha$ not only in  $(-\frac{1}{2\mu^*}, \frac{1}{2\kappa^*})$ but also at $-\frac{1}{2\mu^*}$ and $\frac{1}{2\kappa^*}$.
Although \eq{2.12} looks like a quadratic inequality with respect to $\alpha$ (because of the definitions of $e_0$ and $L$ in \eq{2.6} and \eq{2.4} respectively), we can show that it is linear by expanding $(f_1L_2 + f_2L_1)  e_0 - L_1L_2 \sigma_0$ and seeing that the coefficient of $\alpha^2$ is $0$. Hence, the bound in \eq{2.12} improves or gets worse (depending on the data given) as $\alpha$ tends to $-\frac{1}{2\mu^*}$ or $\frac{1}{2\kappa^*}$ with the optimum value for $\alpha$ occurring at one of the two limits.
The arguments above can be summarized as
\begin{theorem}
The following bound
\beq
	\begin{array}{l}
		\langle\sigma \cdot \varepsilon \rangle - \sigma_0 \cdot \varepsilon_0 \\
		\hspace*{.2in} \geq \min\left\{\frac{1}{f_1f_2}(S_1 - S_2)^{-2}( \varepsilon_0 - \langle S \rangle  \sigma_0) \cdot ((f_1L_{2_*} + f_2L_{1_*})  e_0 - L_{1_*}L_{2_*} \sigma_0) - \frac{a - \det \sigma_0}{\mu^*},\right.\\
		\hspace*{.5in}\left.\frac{1}{f_1f_2}(S_1 - S_2)^{-2}( \varepsilon_0 - \langle S \rangle  \sigma_0) \cdot ((f_1L_{2^*} + f_2L_{1^*})  e_0 - L_{1^*}L_{2^*} \sigma_0) + \frac{a - \det \sigma_0}{\kappa^*}\right\},
	\end{array}
\eeq{2.15} 
where $L_{i_*} = L_i|_{\alpha = -\frac{1}{2\mu^*}}$ and $L_{i^*} = L_i|_{\alpha = \frac{1}{2\kappa^*}},$ $i = 1, 2,$ holds.
\label{theorem 2.1}
\end{theorem}
%%%%%%%%%%%%%%%%%%%%%%%%%%%%%%%%%%%%%%%%%%%%%%%%%%%%%%%%%%%%%%%%%%%%%%%
\section{The attainability condition for \eq{2.15}} \label{The attainability condition for 2.15}
\def\theequation{4.\arabic{equation}}\makeatother
\setcounter{equation}{0}
%%%%%%%%%%%%%%%%%%%%%%%%%%%%%%%%%%%%%%%%%%%%%%%%%%%%%%%%%%%%%%%%%%%%%%%%%%%%%
In this section, we find conditions for the field $\sigma$ to be such that the equality for \eq{2.15} is attained. In other words, fixing $\alpha \in [-\frac{1}{2\kappa^*}, \frac{1}{2\mu^*}],$ we find conditions on $\sigma$ such that
\beq 
	\langle \sigma \cdot L \sigma\rangle  = \langle \hat \sigma \cdot L \hat \sigma \rangle ,
\eeq{3.1} 
where $\hat \sigma$ is a minimizer of the right-hand side in \eq{30} which can be found from \eq{2.7} and \eq{2.10}. 
We have the following lemma.
\begin{lemma}
 	\eq{3.1} holds if, and only if, 
 	\beq
 		L\sigma = L\hat \sigma = \frac{1}{f_1f_2}(L_1 - L_2)^{-1}[f_2\chi_1 L_1(e_0 - L_2\sigma_0) - f_1\chi_2 L_2(e_0 - L_1\sigma_0)].
 	\eeq{3.2}
\label{lemma 3.1}
\end{lemma}

{\it Proof:}
Since $L$ is self-adjoint, the "only if" direction is not hard to see. We thus only prove the "if" direction. Define the functional $J$ by
\beq
	J(\underline \sigma) = \langle \underline \sigma \cdot L \underline \sigma \rangle,
\eeq{} for all vector valued functions $\underline \sigma$ satisfying
\beq
	\langle \underline \sigma \rangle = \sigma_0, \quad \langle L\underline \sigma \rangle = e_0.
\eeq{3.4}   
Since $L$ is semi-positive definite, $J$ is convex. It is not hard to see that for all $t \in [0, 1],$ $\underline \sigma = t\sigma + (1 - t) \hat \sigma$ satisfies \eq{3.4}. Using \eq{3.1} and the fact that $\hat \sigma$ and $\sigma$ are both minimizers of $J$, we have
\beq 
	J(\hat \sigma) \leq J(t\sigma + (1 - t)\hat \sigma) \leq tJ(\sigma) + (1 - t)J(\hat \sigma) = J(\hat \sigma)
\eeq{}  for all $t \in [0, 1].$
It follows that
\beq 
	\frac{d}{dt}J(t\sigma + (1 - t)\hat \sigma) = 0,
\eeq{} 
or equivalently
\beq
	\langle L(t\sigma + (1 - t)\hat \sigma) \cdot (\sigma - \hat \sigma) \rangle = 0.
\eeq{3.7}
Letting $t$ in \eq{3.7} be $1$ and $0$, respectively, gives
\beq 
	\langle L\sigma \cdot (\sigma - \hat \sigma)\rangle  = 0
\eeq{} 
and 
\beq 
	\langle L \hat \sigma \cdot (\sigma - \hat \sigma)\rangle  = 0
\eeq{} 
From the difference of these two equations,
\beq 
	\langle (\sigma - \hat \sigma) \cdot L(\sigma - \hat \sigma)\rangle  = 0,
\eeq{} 
we obtain the first equation of \eq{3.2}.

In order to see the second equation of \eq{3.2}, we calculate $L\hat \sigma$ on each phase. On phase 1,
\beqa
	L\hat \sigma &=& \lambda_1 + L_1\lambda_2\nonum
	&=& 2e_0 - \langle L \rangle\lambda_2 + L_1\lambda_2\nonum
	&=& 2e_0 + f_2(L_1 - L_2)\lambda_2\nonum
	&=& 2e_0 + \frac{2}{f_1}(L_1 - L_2)^{-1}[(f_2L_1 + f_1L_2)e_0 - L_1L_2\sigma_0]\nonum
	&=& \frac{2}{f_1}(L_1 - L_2)^{-1}[f_1(L_1 - L_2)e_0 + (f_2L_1 + f_1L_2)e_0 - L_1L_2\sigma_0]\nonum
	&=& \frac{2}{f_1}(L_1 - L_2)^{-1}L_1(e_0 - L_2\sigma_0).   
\eeqa{}
Similarly, on phase 2,
\beq
	L\hat \sigma = - \frac{2}{f_2}(L_1 - L_2)^{-1}L_2(e_0 - L_1\sigma_0).
\eeq{} 
The lemma follows.
\qed

Lemma \ref{lemma 3.1} implies that if bound \eq{2.15} is attained then there are two constant vectors $D$ and $E$ such that
\beq
	L \sigma = L \hat \sigma = \chi_1 D + \chi_2 E.
\eeq{3.11}
We next show that \eq{3.11} is the attainability condition for \eq{2.15}, that we are looking for in this section. In fact, denoting by $\langle \sigma \rangle_1$ and $\langle \sigma \rangle_2$ the average of $\sigma$ on phase $1$ and phase $2$, respectively, as in the previous section, we have
\beq
	D = L_1 \langle \sigma \rangle_1,\quad E = L_2 \langle \sigma \rangle_2.
\eeq{3.12}

Plugging \eq{2.11} into \eq{3.12} and then the resulting values of $D$ and $E$ into \eq{3.11}, we obtain \eq{3.2}.

We have proved the theorem.
\begin{theorem} 
	The bound in \eq{2.15} becomes equality if, and only if,
		\beq
			L_* \sigma = \chi_1D + \chi_2 E
		\eeq{} where $D$ and $E$ are two constant vectors and $L_*$ is either $L|_{\alpha = -\frac{1}{2\mu^*}}$ or $L|_{\alpha = \frac{1}{2\kappa^*}}$. Moreover, 
	\begin{enumerate}
		\item if $\mu^* = \mu_1\ne\mu_2$ then
			\beq
				\begin{array}{l}
					\langle\sigma \cdot \varepsilon \rangle - \sigma_0 \cdot \varepsilon_0 \\
		\hspace*{.2in} = \frac{1}{f_1f_2}(S_1 - S_2)^{-2}( \varepsilon_0 - \langle S \rangle  \sigma_0) \cdot ((f_1L_{2_*} + f_2L_{1_*})  e_0 - L_{1_*}L_{2_*} \sigma_0) - \frac{a - \det \sigma_0}{\mu^*},
				\end{array}
			\eeq{} is equivalent to 
			\beq
				L|_{\alpha = -\frac{1}{2\mu_1}}\sigma = \chi_1 D_* + \chi_2 E
			\eeq{} with $D_* = (d_1, 0, 0)$ and this holds if and only if the field $\sigma$ is constant on phase 2, and has constant first component (bulk part) on phase 1;
		\item if $\kappa^* = \kappa_1\ne\kappa_2$ then
		\beq
	\begin{array}{l}
		\langle\sigma \cdot \varepsilon \rangle - \sigma_0 \cdot \varepsilon_0 \\
		\hspace*{.2in} = \frac{1}{f_1f_2}(S_1 - S_2)^{-2}( \varepsilon_0 - \langle S \rangle  \sigma_0) \cdot ((f_1L_{2^*} + f_2L_{1^*})  e_0 - L_{1^*}L_{2^*} \sigma_0) + \frac{a - \det \sigma_0}{\kappa^*}
	\end{array}
\eeq{} is equivalent to
	\beq
		L|_{\alpha = \frac{1}{2\kappa_1}}\sigma = \chi_1 D_* + \chi_2 E
	\eeq{} with $D_* = (0, d_2, d_3)$ and this holds if and only if the field $\sigma$ is constant on phase 2, and has constant second and third components  (shear part) on phase 1;.
	\end{enumerate}
\label{theorem 3.2}
\end{theorem}

%%%%%%%%%%%%%%%%%%%%%%%%%%%%%%%%%%%%%%%%%%%%%%%%%%%
\section{The method of splitting}\label{The method of splitting}
\def\theequation{5.\arabic{equation}}\makeatother
\setcounter{equation}{0}
%%%%%%%%%%%%%%%%%%%%%%%%%%%%%%%%%%%%%%%%%%%%%%%%%%%
We introduce here another approach, not based on variational principles, using which we can deduce the previous bounds. The main idea of the method is to split the domain
into its phases, to correlate averages over each phase, and then
to project the information to obtain the desired bound. Moreover, this method allows us to add one more datum: the null-Lagrangian $b = \langle \det \nabla u\rangle.$ In other words, the known quantities are
\beq
	E = \langle \sigma \cdot \nabla u \rangle, \quad \sigma_0 = \langle \sigma \rangle, \quad \langle \nabla u \rangle, \quad a = \langle \det \sigma \rangle, \quad 	b = \langle \det \nabla u\rangle
\eeq{}  

Consider the four dimensional space of $2 \times 2$ matrices with the basis
\beq
	B = \left\{
			\frac{1}{\sqrt{2}}\left[
								\begin{array}{cc}
									0&1\\-1&0
								\end{array}
							  \right],
			\frac{1}{\sqrt{2}}\left[
								\begin{array}{cc}
									1&0\\0&1
								\end{array}
							  \right],
			\frac{1}{\sqrt{2}}\left[
								\begin{array}{cc}
									1&0\\0&-1
								\end{array}
							  \right],
			\frac{1}{\sqrt{2}}\left[
								\begin{array}{cc}
									0&1\\1&0
								\end{array}
							  \right]
		  \right\},
\eeq{} so that 
\beq
	\sigma = (0, \sigma_1, \sigma_2, \sigma_3), \quad \nabla u = (F_0, \varepsilon_1, \varepsilon_2, \varepsilon_3) 
\eeq{} and
\beq
	a = \frac{1}{2}\langle \sigma_1^2 - \sigma_2^2 - \sigma_3^2\rangle, \quad c \equiv b - \frac{1}{2}\langle F_0\rangle^2 \geq b - \frac{1}{2}\langle F_0^2 \rangle = \frac{1}{2}\langle \varepsilon_1^2 - \varepsilon_2^2 - \varepsilon_3^2\rangle.
\eeq{6.5}
Note that the inequality in \eq{6.5} is attained if, and only if, $F_0$ is constant everywhere.

Defining (b for bulk and s for shear)
\beq
	\begin{array}{rcl}
		E_{1b} &=& \langle \chi_1 \sigma_1\varepsilon_1 \rangle = 2\kappa_1\langle \chi_1\varepsilon_1^2\rangle, \\
		E_{2b} &=& \langle \chi_2 \sigma_1\varepsilon_1 \rangle = 2\kappa_2\langle \chi_2\varepsilon_1^2\rangle, \\
		E_{1s} &=& \langle\chi_1(\sigma_2\varepsilon_2 + \sigma_3\varepsilon_3)\rangle = 2\mu_1\langle\chi_1(\varepsilon_2^2 + \varepsilon_3^2)\rangle,\\
		E_{2s} &=& \langle\chi_2(\sigma_2\varepsilon_2 + \sigma_3\varepsilon_3)\rangle = 2\mu_2\langle\chi_2(\varepsilon_2^2 + \varepsilon_3^2)\rangle,	
	\end{array}
\eeq{}
we have
\beq	
	E = E_{1b} + E_{1s} + E_{2b} + E_{2s},
\eeq{6.7}
\beq
	a = \kappa_1E_{1b} + \kappa_2E_{2b} - \mu_1E_{1s} - \mu_2E_{2s},
\eeq{6.8}
\beq
	c \geq \frac{E_{1b}}{4\kappa_2} + \frac{E_{2b}}{4\kappa_2} - \frac{E_{1s}}{4\mu_1} - \frac{E_{2s}}{4\mu_2}.
\eeq{6.9}

On the other hand, it follows from the identities
\beqa{}
	\langle \varepsilon_i \rangle & = & \langle \chi_1 \varepsilon_i\rangle + \langle \chi_2\varepsilon_i\rangle,   \quad i = 1, 2, 3, \nonum
	\langle \sigma_1 \rangle & = & 2\kappa_1 \langle \chi_1\varepsilon_1\rangle + 2\kappa_2\langle \chi_2 \varepsilon_2\rangle, \quad 
	\langle \sigma_j \rangle = 2\mu_1 \langle \chi_1\varepsilon_j\rangle + 2\mu_2\langle \chi_2 \varepsilon_j\rangle, \quad j = 2, 3,
\eeqa{}
that 
\beqa
	\langle \chi_1\varepsilon_1\rangle & = & \frac{1}{2(\kappa_2 - \kappa_1)}(2\kappa_2\langle\varepsilon_1\rangle - \langle \sigma_1\rangle), \quad \langle \chi_2\varepsilon_1\rangle = \frac{1}{2(\kappa_1 - \kappa_2)}(2\kappa_1\langle\varepsilon_1\rangle - \langle \sigma_1
	\rangle), \nonum
	\langle\chi_1 \varepsilon_j\rangle & = & \frac{1}{2(\mu_2 - \mu_1)}(2\mu_2\langle \varepsilon_j\rangle - \langle \sigma_j\rangle), \quad \langle\chi_1 \varepsilon_j\rangle = \frac{1}{2(\mu_2 - \mu_1)}(2\mu_2\langle \varepsilon_j\rangle - \langle \sigma_j\rangle), \quad j = 2, 3.
\eeqa{6.12} 
So, these quantities are known.

\begin{figure}[htp]
\centering
	\begin{picture}(100, 130)
			\put(-50,0){\includegraphics[scale=.30]{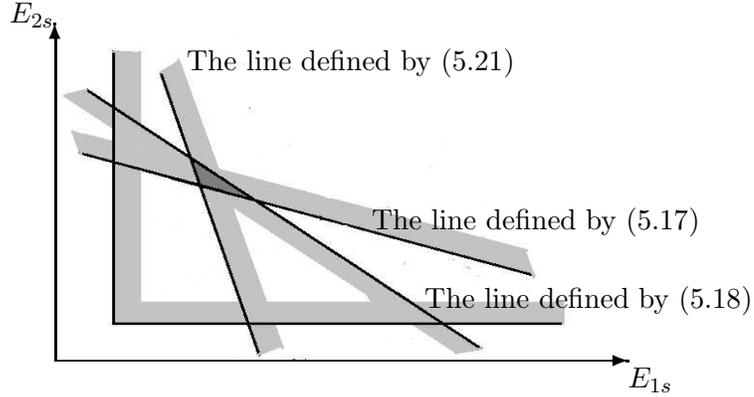}}
			\thicklines
			\put(166,0){\vector(1,0){1}}
			\put(-50,127){\vector(0,1){1}}
			\put(-67, 128){$E_{2s}$}
			\put(167,-10){$E_{1s}$}
			\put(70,50){The line defined by \eq{6.19}}
			\put(90,20){The line defined by \eq{6.20}}
			\put(0,110){The line defined by \eq{6.22}}
		\end{picture}
\caption{The region of possible $(E_{1s}, E_{2s}),$ assuming that $\kappa_1 > \kappa_2$.}
\label{figure 3}
\end{figure}

\begin{figure}[htp]
\centering
	\begin{picture}(100, 130)
			\put(-50,0){\includegraphics[scale=.30]{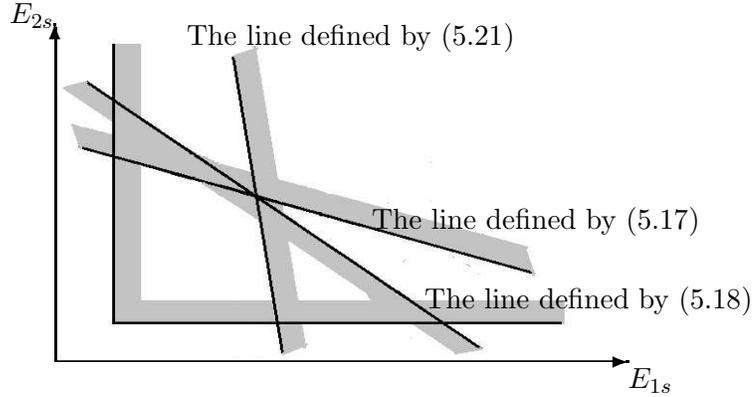}}
			\thicklines
			\put(166,0){\vector(1,0){1}}
			\put(-50,127){\vector(0,1){1}}
			\put(-67, 128){$E_{2s}$}
			\put(167,-10){$E_{1s}$}
			\put(70,50){The line defined by \eq{6.19}}
			\put(90,20){The line defined by \eq{6.20}}
			\put(0,120){The line defined by \eq{6.22}}
		\end{picture}
\caption{The feasible region shrinking to a point, assuming that $\kappa_1 > \kappa_2$.}
\label{figure 4}
\end{figure}

Note that for $i = 1, 2, 3,$
\beq
	\langle \chi_1\varepsilon_i^2\rangle - \frac{1}{f_1}\langle\chi_1\varepsilon_i\rangle^2 =  \langle(\chi_1 \varepsilon_i - \frac{\chi_1}{f_1}\langle \chi_1\varepsilon_k\rangle)^2\rangle \geq 0,
\eeq{6.13} 
with equality when $\varepsilon_i$ is constant on phase $1$. Similarly,
\beq
	\langle \chi_2\varepsilon_i^2\rangle - \frac{1}{f_2}\langle\chi_2\varepsilon_i\rangle^2 \geq 0, \quad i = 1, 2, 3,
\eeq{6.14} 
with equality when $\varepsilon_i$ is constant on phase 2. Therefore, defining the known quantities
\beq
	\begin{array}{rcl}
		A_{1b} &=& 2\kappa_1\langle\chi_1\varepsilon_1\rangle^2,\\
		A_{2b} &=& 2\kappa_2\langle\chi_2\varepsilon_1\rangle^2,\\
		A_{1s} &=& 2\mu_1(\langle \chi_1 \varepsilon_2\rangle^2 + \langle \chi_1\varepsilon_3 \rangle^2),\\
		A_{2s} &=& 2\mu_2(\langle \chi_2 \varepsilon_2\rangle^2 + \langle \chi_2\varepsilon_3 \rangle^2),
	\end{array}
\eeq{6.15}
we have from \eq{6.13} and \eq{6.14} 
\beq
	E_{1b} \geq \frac{A_{1b}}{f_1}, \quad E_{2b} \geq \frac{A_{2b}}{f_2},
\eeq{6.16}
\beq
	E_{1s} \geq \frac{A_{1s}}{f_1}, \quad E_{2s} \geq \frac{A_{2s}}{f_2}.
\eeq{6.17}
Also, $E_{1b}$ and $E_{2b}$ can be eliminated using \eq{6.7} and \eq{6.8}
\beq
	\begin{array}{rcl}
		E_{1b} &=& \frac{1}{\kappa_1 - \kappa_2}(a - \kappa_2 E + (\kappa_2 + \mu_1)E_{1s} + (\kappa_2 + \mu_2)E_{2s}),\\
		E_{2b} &=& \frac{1}{\kappa_2 - \kappa_1}(a - \kappa_1 E + (\kappa_1 + \mu_1)E_{1s} + (\kappa_1 + \mu_2)E_{2s}).
	\end{array}
\eeq{}
	So, \eq{6.16} gets replaced by
\beq
	\frac{1}{\kappa_1 - \kappa_2}(a - \kappa_2 E + (\kappa_2 + \mu_1)E_{1s} + (\kappa_2 + \mu_2)E_{2s}) \geq \frac{A_{1b}}{f_1},
\eeq{6.19} and
\beq
\frac{1}{\kappa_2 - \kappa_1}(a - \kappa_1 E + (\kappa_1 + \mu_1)E_{1s} + (\kappa_1 + \mu_2)E_{2s}) \geq \frac{A_{2b}}{f_2}
\eeq{6.20} and the inequality in \eq{6.9} gets replaced by
\beq
	\begin{array}{rcl}
		c &\geq& \frac{1}{4\kappa_1(\kappa_1 - \kappa_2)}(a - \kappa_2 E + (\kappa_2 + \mu_1)E_{1s} + (\kappa_2 + \mu_2)E_{2s})\\
		&&\hspace*{.24in} + \frac{1}{4\kappa_2(\kappa_2 - \kappa_1)} (a - \kappa_1 E + (\kappa_1 + \mu_1)E_{1s} + (\kappa_1 + \mu_2)E_{2s})%\\
		%&& \hspace*{.24in}
		-\frac{E_{1s}}{4\mu_1} - \frac{E_{2s}}{4\mu_2}.
	\end{array}
\eeq{6.21}
Now we have
\beqa
	\frac{1}{\kappa_1(\kappa_1 - \kappa_2)} + \frac{1}{\kappa_2(\kappa_2 - \kappa_1)} &=& -\frac{1}{\kappa_1\kappa_2},\nonum
	\frac{-\kappa_2}{\kappa_1(\kappa_1 - \kappa_2)} - \frac{\kappa_1}{\kappa_2(\kappa_2 - \kappa_1)} &=& \frac{\kappa_1 + \kappa_2}{\kappa_1\kappa_2},\nonum
	\frac{\kappa_2 + \mu_1}{\kappa_1(\kappa_1 - \kappa_2)} + \frac{\kappa_1 + \mu_1}{\kappa_2(\kappa_2 - \kappa_1)} - \frac{1}{\mu_1} &=& \frac{\mu_1\kappa_2(\kappa_2 + \mu_1) - \mu_1 \kappa_1(\kappa_1 + \mu_1) - \kappa_1\kappa_2(\kappa_1 - \kappa_2)}{\kappa_1\kappa_2\mu_1(\kappa_1 - \kappa_2)}\nonum
	&=& \frac{-\mu_1(\kappa_1 + \kappa_2) - \mu_1^2 - \kappa_1\kappa_2}{\kappa_1\kappa_2\mu_1} \nonum
	&=& -\frac{(\kappa_1 + \mu_1)(\kappa_2 + \mu_1)}{\kappa_1\kappa_2\mu_1}, \nonum
	\frac{\kappa_1 + \mu_2}{\kappa_1(\kappa_1 - \kappa_2)} + \frac{\kappa_1 + \mu_2}{\kappa_2(\kappa_2 - \kappa_1)} - \frac{1}{\mu_2} 
	&=& -\frac{(\kappa_1 + \mu_2)(\kappa_2 + \mu_2)}{\kappa_1\kappa_2\mu_2}.
\eeqa{}
So, \eq{6.21} becomes
\beq
a - E(\kappa_1 + \kappa_2) + \frac{E_{1s}(\kappa_1 + \mu_1)(\kappa_2 + \mu_1)}{\mu_1} + \frac{E_{2s}(\kappa_1 + \mu_2)(\kappa_2 + \mu_2)}{\mu_2} \geq -4\kappa_1\kappa_2 c.
\eeq{6.22}
So, \eq{6.17}, \eq{6.19}, \eq{6.20} and \eq{6.22} give us $5$ inequalities on the pair $(E_{1s}, E_{2s}).$ To project out the information about the unknowns $(E_{1s}, E_{2s})$ we observe
that the volume fraction $f_1 = 1 - f_2$ must be such that these inequalities define a feasible region in the $(E_{1s}, E_{2s})$ plane (see e.g. Figure \ref{figure 3}).
Note that when either $f_1$ or $f_2$ goes to $0$, \eq{6.17} cannot be satisfied. In this case, the feasible region is empty. Thus, in the generic case, at a bound on $f_1$, i.e. at the limiting value of $f_1$, the feasible region shrinks to a point. In other words, in the generic case, $3$ of the inequalities will be satisfied as equalities and the remaining $2$ as inequalities: the picture looks like, e.g. Figure \ref{figure 4}.
Note that the 3 inequalities which are satisfied as equalities must have "outward normals" $n_1, n_2, n_3$ such that
\beq
	\alpha_1n_1 + \alpha_2n_2 + \alpha_3n_3 = 0
\eeq{6.23} for some $\alpha_1, \alpha_2, \alpha_3 \geq 0.$
The previous approach in Section \ref{A bound obtained by using the null-Lagragian} only took the four inequalities in \eq{6.17}, \eq{6.19} and \eq{6.20} into account. The limiting value of $f_1$ corresponds to the case where $3$ of these were satisfied as equalities and the fourth as an inequality. Thus, one sees immediately that the attainability condition is where the field is constant in one phase and either the bulk part or shear part is constant in the other phase.
With the additional inequality \eq{6.22}, there will be $10$ cases to consider (number of ways to pick $3$ equations from $5$), $4$ considered above and $6$ new.
However, not all of these 10 cases satisfy \eq{6.23}; e.g. if $\kappa_1 > \kappa_2$, the triplet of equations that correspond to \eq{6.17} and \eq{6.22} do not satisfy \eq{6.23} since if $E_{1s}$ and $E_{2s}$ are both large and positive, all three equations are satisfied; i.e. the feasible region given by
 the three equations cannot be a point, but is instead an open region. The question arises as to which of these 6 new cases satisfy \eq{6.23}.

It is obvious that the outward normals of the regions that satisfy the two inequalities in \eq{6.17} are, respectively, 
\beq \nu_1 = (-1, 0)\eeq{} and \beq \nu_2 = (0, -1).\eeq{}  
Without loss of generality, by relabeling the phases if necessary, assume that $\kappa_1 > \kappa_2.$ In this case, an outward normal of the region that satisfies \eq{6.19} is 
\beq \nu_3 = (-\kappa_2 - \mu_1, -\kappa_2 - \mu_2),\eeq{} 
an outward normal of the region that satisfies \eq{6.20} is 
\beq \nu_4 = (\kappa_1 + \mu_1, \kappa_1 + \mu_2), \eeq{} 
and, finally, an outward normal of the region that satisfies \eq{6.22} is 
\beq\nu_5 = \left(-\frac{(\kappa_1 + \mu_1)(\kappa_2 + \mu_1)}{\mu_1}, -\frac{(\kappa_1 + \mu_2)(\kappa_2 + \mu_2)}{\mu_2}\right).\eeq{} 
We now consider each triplet of normals that include $\nu_5$ to see when they satisfy \eq{6.23}. As said, $(\nu_1, \nu_2, \nu_5)$ does not satisfy \eq{6.23}. If $(\alpha, \beta)$ solves 
\beq \nu_5 = \alpha \nu_1 + \beta \nu_3,\eeq{} 
then
\beq
	\beta = \frac{\kappa_1 + \mu_2}{\mu_2} > 0. 
\eeq{} Thus, the triplet $(\nu_1, \nu_3, \nu_5)$
 never satisfies \eq{6.23}. Similarly, neither does the triplet $(\nu_2, \nu_3, \nu_5)$. Next, solving the equation
 \beq \nu_5 = \alpha \nu_1 + \beta \nu_4,\eeq{}
gives
\beq
	\alpha = -\frac{\kappa_2(\mu_1 - \mu_2)(\kappa_1 + \mu_1)}{\mu_1\mu_2}, \quad \beta = -\frac{\kappa_2 + \mu_2}{\mu_2}.
\eeq{6.33} So, the triplet $(\nu_1, \nu_4, \nu_5)$ satisfies \eq{6.23} if, and only if, $\mu_1 > \mu_2.$  By a similar technique, solving 
\beq \nu_5 = \alpha \nu_2 + \beta \nu_4,\eeq{}
we have
\beq
	\alpha = -\frac{\kappa_2(\kappa_1 + \mu_2)(\mu_2 - \mu_1)}{\mu_1\mu_2}, \quad \beta = -\frac{\kappa_2 + \mu_1}{\mu_1}.
\eeq{} Therefore, the triplet $(\nu_2, \nu_4, \nu_5)$ satisfies \eq{6.23} if, and only if, $\mu_1 < \mu_2.$ Finally, we consider the triplet $(\nu_3, \nu_4, \nu_5)$. The unique solution of
\beq
	\nu_5 = \alpha \nu_3 + \beta \nu_4
\eeq{} is
\beq
	\alpha = -\frac{\kappa_2(\kappa_1 + \mu_1)(\kappa_1 + \mu_2)}{(\kappa_1 - \kappa_2)\mu_1\mu_2} < 0, \quad \beta = -\frac{\kappa_1(\kappa_2 + \mu_1)(\kappa_2 + \mu_2)}{(\kappa_1 - \kappa_2)\mu_1\mu_2} < 0. \eeq{6.37} Hence, this triplet always satisfies \eq{6.23}.

We next find the desired bounds. We can do so by finding an appropriate linear combination of the equations in each triplet involving \eq{6.22}
that satisfies \eq{6.23} to obtain
\begin{theorem} 
Assuming $\kappa_1 > \kappa_2$, and recalling the definitions \eq{6.15} and \eq{6.12} of $A_{1b}$, $A_{2b}$, $A_{1s}$, and $A_{2s}$ 
one has the bound
\beq
4 c \kappa_1 \kappa_2 \geq -\frac{\kappa_1\kappa_2(a + E(\mu_1 + \mu_2))}{\mu_1\mu_2} +\frac{A_{1b} \kappa_2 (\kappa_1+\mu_1) (\kappa_1+\mu_2)}{f_1 \mu_1 \mu_2}+\frac{A_{2b} \kappa_1 (\kappa_2+\mu_1) (\kappa_2+\mu_2)}{f_2 \mu_1 \mu_2},
\eeq{6.38}
with equality when $F_0$ (the component of the antisymmetric part of $\nabla u$)
is constant everywhere, and the bulk component $\varepsilon_1$ is constant in phase 1 and constant in phase2.
If additionally $\mu_1 > \mu_2$ then one has
\beq
	4 c \kappa_1 \kappa_2 \geq \frac{\kappa_2 (a+E (-\kappa_1+\mu_2))}{\mu_2} +\frac{A_{1s} \kappa_2 (\kappa_1+\mu_1) (\mu_1-\mu_2)}{f_1 \mu_1 \mu_2}+\frac{A_{2b} (\kappa_1-\kappa_2) (\kappa_2+\mu_2)}{f_2 \mu_2},
\eeq{6.39}
with equality when $F_0$ is constant everywhere, the shear components $\varepsilon_2$ and $\varepsilon_3$ are constant in phase 1, and the bulk component $\varepsilon_1$ is constant in phase 2.
Alternatively if $\mu_1 < \mu_2$ (and $\kappa_1 > \kappa_2$) 
then one has
\beq
	4 c \kappa_1 \kappa_2 \geq \frac{\kappa_2 (a+E (-\kappa_1+\mu_1))}{\mu_1} +\frac{A_{2s} \kappa_2 (\mu_2-\mu_1) (\kappa_1+\mu_2)}{f_1 \mu_1 \mu_2} + \frac{A_{2b} (\kappa_1-\kappa_2) (\kappa_2+\mu_1)}{f_2 \mu_1},
\eeq{6.40}
with equality when $F_0$ is constant everywhere, and both bulk and shear components $\varepsilon_1$, $\varepsilon_2$ and $\varepsilon_3$ are constant in phase 2.
\label{theorem 4.1}
\end{theorem}

{\it Proof:}
As suggested by \eq{6.37}, we can add \eq{6.22} to \eq{6.19}, multiplied by $\ds\frac{\kappa_2 (\kappa_1+\mu_1) (\kappa_1+ \mu_2)}{\mu_1 \mu_2},$ and \eq{6.20}, multiplied by $\ds \frac{\kappa_1 (\kappa_2+ \mu_1) (\kappa_2+ \mu_2)}{\mu_1 \mu_2},$ to eliminate both $E_{1s}$ and $E_{2s}$. Doing so, we obtain \eq{6.38}. Similarly in the case that $\mu_1 > \mu_2$, we can use the first inequality in \eq{6.17}, \eq{6.20} and \eq{6.22} to deduce
\eq{6.39}. Finally, if $\mu_1 < \mu_2$ then the second inequality in \eq{6.17}, \eq{6.20} and \eq{6.22} yield \eq{6.40}. The attainability conditions follow directly from the conditions under which the inequalities
in \eq{6.5} and \eq{6.13} are satisfied as equalities. 
\qed

It is an open question as to whether the bounds of Theorem \ref{theorem 4.1} could be obtained using the translation method. Although we suspect they could, the method of splitting 
has the advantage of immediately providing the attainability conditions. 

\begin{remark}
The final set of bounds are the intersection of the inequalities provided by Theorem \ref{theorem 2.1} and Theorem \ref{theorem 4.1}. Each of these inequalities can be multiplied by $f_1f_2$ to yield quadratic inequalities
in $f_1=1-f_2$ which may be easily solved to give the maximum interval of $f_1$ compatible with all the inequalities. If one is interested in bounds on (average stress, average strain) pairs in composites
then one should make the substitutions \eq{comp} in these inequalities. It remains to be investigated whether the resulting bounds provide a complete characterization of the possible 
(average stress, average strain) pairs in composites of two phases mixed in a given proportion, or whether there are some missing bounds. 
\end{remark} 

\section*{Acknowledgements}
The authors thank Andre Zaoui for his many contributions to micromechanics, and
are grateful for support from the National Science Foundation through grant DMS-0707978.
\ifx \bblindex \undefined \def \bblindex #1{} \fi\ifx \bblindex \undefined \def
  \bblindex #1{} \fi

%%%%%%%%%%%%%%%%%%%%%%%%%%%%%%%%%%%%%%%%%%%%%%%%  

%\bibliography{/u/ma/milton/tcbook,/u/ma/milton/newref} 
%\bibliographystyle{unsrt}

\end{document}